\documentclass[11pt, a4paper]{amsart}
\usepackage{amscd, amssymb, pb-diagram} 
\input{prepictex.tex}
\input{pictex.tex}
\input{postpictex.tex}

\def\T{\mathbb T}
\def\Z{\mathbb Z}

\def\C{\mathbb C}

\def\H{{\mathcal H}}
\def\K{{\mathcal K}}

\def\O{{\mathcal O}}
\def\Hbad{H^{\text{{\rm fin}}}_{\infty}}
\def\Kbad{K^{\text{{\rm fin}}}_{\infty}}
\def\Hibad{(H_i)^{\text{{\rm fin}}}_{\infty}}

\def\Obad{\Omega (v)^{\text{{\rm fin}}}_{\infty}}

\def\Aut{\operatorname{Aut}}
\def\coker{\operatorname{coker}}

\newtheorem{thm}{Theorem}[section]
\newtheorem{cor}[thm]{Corollary}
\newtheorem{lemma}[thm]{Lemma}
\newtheorem{prop}[thm]{Proposition}

\theoremstyle{definition}

\theoremstyle{remark}
\newtheorem{remark}[thm]{Remark}
\newtheorem{example}[thm]{Example}

\begin{document} 

\title[Ideal Structure of Graph Algebras]
      {\boldmath{The Ideal Structure of the\\
$C^*$-Algebras of Infinite Graphs}} 

\author[Bates]{Teresa Bates} 
\address{School of Mathematics\\ 
The University of New South Wales\\ 
Sydney, NSW 2052\\ 
Australia} 
\email{teresa@maths.unsw.edu.au} 

\author[Hong]{Jeong Hee Hong}
\address{Applied Mathematics\\ 
Korea Maritime University\\ 
Pusan 606--791\\ 
South Korea} 
\email{hongjh@hanara.kmaritime.ac.kr} 

\author[Raeburn]{Iain Raeburn}
\address{Mathematics\\
The University of Newcastle\\
Callaghan, NSW 2308\\
Australia}
\email{iain@maths.newcastle.edu.au} 

\author[Szyma\'{n}ski]{Wojciech Szyma\'{n}ski}
\address{Mathematics\\
The University of Newcastle\\
Callaghan, NSW 2308\\
Australia}
\email{wojciech@frey.newcastle.edu.au} 

\thanks{This research was partially supported by the 
Korea Maritime University Research Foundation, the University of 
Newcastle, and the Australian Research Council.} 

\date{11 September 2001} 

\maketitle

\begin{abstract}
We classify the gauge-invariant ideals in the $C^*$-algebras of 
infinite directed graphs, and describe the quotients as graph 
algebras. We then use these results to identify the gauge-invariant 
primitive ideals in terms of the structural properties of the graph, 
and describe the $K$-theory of the $C^*$-algebras of arbitrary 
infinite graphs. 
\end{abstract}

\section{Introduction} 

There has recently been a great deal of interest in generalisations 
of the Cuntz-Krieger algebras associated to infinite directed graphs 
\cite{kprr,flr} and infinite matrices \cite{el1}. The basic theorems 
of Cuntz and Krieger \cite{ck,c} on uniqueness and ideal structure 
have elegant extensions to the $C^*$-algebras of the row-finite graphs 
in which each vertex emits only finitely many edges 
\cite{kprr,kpr,kpw,bprs}. Various authors have investigated 
the $C^*$-algebras of arbitrary infinite graphs from different points of 
view, obtaining satisfactory versions of the uniqueness theorems
\cite{flr,rs,sz4}. However, these articles do not provide a complete 
description of the ideal structure of graph algebras, as 
is given in \cite{hr} for the Cuntz-Krieger algebras of 
finite matrices. Indeed, even for row-finite graphs 
the ideal structure has only been well-understood when the graph satisfies 
the Condition (K) introduced in \cite{kprr} (see \cite{bprs}). 

The analysis in \cite{hr} shows that to understand the ideal structure 
of graph algebras we first need to describe the gauge-invariant 
ideals. The main purpose of this paper is to provide such a 
description for arbitrary infinite graphs. We give a complete 
list of the gauge-invariant ideals of $C^*(E)$ for an arbitrary 
infinite graph $E$ (Theorem \ref{gaugeinvariant}), and 
then use it to identify all the gauge-invariant primitive ideals 
(Theorem \ref{primitivegaugeinvariant}). When the 
graph satisfies Condition (K) all ideals are gauge-invariant 
and our results give their complete classification. 

The key tool in our approach is a realisation of the quotient 
$C^*(E)/J$ by a gauge-invariant ideal as the graph algebra of a quotient graph 
(Proposition \ref{quotientalgebra}). This result is  of considerable interest
in its own right, because we are able to explicitly describe the quotient
graph. As a further application, we show how to extend the description  of 
$K_*(C^*(E))$ obtained in \cite{rs} for row-finite
$E$ to arbitrary infinite graphs 
(Theorem \ref{Kgroups}). 

There are several reasons for the current interest in graph algebras 
apart from the elegance of their theory. First, they provide 
good test problems in the general theories of groupoid algebras 
\cite{kprr,kpr,pat}, Cuntz-Pimsner algebras \cite{pin,sch,kpw,flr,fmr}, and 
partial crossed products \cite{el1,elq}. Second, the simple graph 
algebras provide a rich family of accessible models for  
purely infinite simple $C^*$-algebras. Indeed, Szyma\'{n}ski has 
shown in \cite{sz3} that every stable, purely infinite, simple and classifiable 
$C^*$-algebra with $K_1$ torsion-free  can be realised 
as a graph algebra. Although there is some debate about what `purely infinite' 
should mean for non-simple $C^*$-algebras \cite{kr}, there is already
considerable  interest in their classification, and it is likely that the
non-simple  graph algebras will again provide an important family of models. 

We have had the main results of the present paper for some time 
(see \cite{h}), and had wanted to include them in a complete 
analysis of the ideal structure of infinite graph algebras. However, 
we have received many enquires about this work, and in response have 
decided to publish it in stages. There is some overlap between the 
present article and the work of Drinen and Tomforde \cite{dt1}, 
who describe the primitive ideal space of the $C^*$-algebras of 
graphs satisfying Condition (K) by reducing to the row-finite 
case.  Our methods are quite different from theirs: we  work directly with
quotients of graph algebras rather than algebras 
Morita equivalent to them. In the sequel,  we use
these techniques to obtain a complete generalisation of  the program of
\cite{hr} to the
$C^*$-algebras of arbitrary infinite graphs. 

\section{Preliminaries} 

Let $E=(E^0,E^1,r,s)$ be a (countable) directed graph, consisting of a set
$E^0$ of vertices, a set $E^1$ of edges, and range and source maps $r,s:E^1\to
E^0$. A Cuntz-Krieger
$E$-family consists of mutually orthogonal projections
$\{P_v:v\in E^0\}$ and partial isometries $\{S_e:e\in
E^1\}$ with mutually orthogonal ranges satisfying
\begin{description}
\item[(G1)] $S_e^*S_e=P_{r(e)}$,
\item[(G2)] $S_eS_e^*\leq P_{s(e)}$, and
\item[(G3)] $P_v=\sum_{s(e)=v}S_eS_e^*$  if $s^{-1}(v)$ is
finite and non-empty.
\end{description}
The $C^*$-algebra $C^*(E)$ of  $E$ is  the universal
$C^*$-algebra generated by a Cuntz-Krieger $E$-family
$\{s_e,p_v\}$. If $\{S_e,P_v\}$ is a Cuntz-Krieger $E$-family, we
denote by $\pi_{S,P}$ the representation of $C^*(E)$  such that
$\pi_{S,P}(p_v)=P_v$ and $\pi_{S,P}(s_e)=S_e$.

We denote by $\gamma:\T\rightarrow\Aut C^*(E)$ the gauge 
action, which is characterised on generators by $\gamma_z(p_v)=p_v$ and
$\gamma_z (s_e)=zs_e$ for $v\in E^0$, $e\in E^1$, $z\in\T$. Existence of 
the gauge action is equivalent to universality in the definition 
of $C^*(E)$, as the following gauge-invariant uniqueness theorem 
shows. This result was proved for finite graphs in
\cite[Theorem~2.3]{hr}, for row-finite graphs in 
\cite[Theorem 2.1]{bprs}, and generalised in \cite[Theorem 2.7]{rs} 
to the Cuntz-Krieger algebras of infinite matrices and in
\cite[Theorem~4.1]{fmr} to Cuntz-Pimsner algebras. Unfortunately, the  existing
versions do not  cover all infinite graphs with sources or sinks. 

\begin{thm}\label{uniqueness}  
Let $E$ be an arbitrary directed graph, let $\{S_e,P_v\}\subset B(\H_E)$ be a
Cuntz-Krieger 
$E$-family, and let $\pi=\pi_{S,P}$ be the representation of $C^*(E)$ such 
that $\pi(s_e)=S_e$ and $\pi(p_v)=P_v$. Suppose that each $P_v$ is non-zero,
and that there is a strongly continuous action $\beta$ of $\T$ on
$C^*(S_e,P_v)$ such that $\beta_z \circ \pi=\pi\circ \gamma_z$ for $z\in
\T$. Then $\pi$ is faithful. 
\end{thm}
\begin{proof} 
If $E$ is an infinite directed graph without sinks or sources (that is, 
each vertex emits and receives some edges), then $C^*(E)$ is naturally 
isomorphic to the Cuntz-Krieger algebra of a suitable infinite matrix 
\cite[Theorem~10]{flr}. Thus \cite[Theorem 2.7]{rs} applies to such graphs. 

To extend the theorem to graphs with sinks, it suffices to 
add tails as in \cite{bprs}. Indeed, let $F$ be the graph obtained by 
adding a tail (with extra vertices $\{v_i:i=1,2,\ldots\}$) 
to a sink $w$ of $E$ as in \cite[\S1]{bprs}, let
$\H_T=\bigoplus_{i=1}^\infty \H_i$ be the direct sum of copies $\H_i$ of
$P_w\H_E$, let 
$\{T_e,Q_v\}$ be the Cuntz-Krieger $F$-family on $\H_F=\H_E\oplus\H_T$ 
obtained by extending the Cuntz-Krieger $E$-family $\{S_e,P_v\}$ as in 
\cite[Lemma 1.2]{bprs}, and let $U:\T\rightarrow{U}(\H_E)$ 
be a unitary representation such that $(\pi_{S,P},U)$ is covariant 
for the gauge action on $C^*(E)$. Then there is a unitary representation 
$V:\T\rightarrow{U}(\H_F)$ such that $(\pi_{T,Q},V)$ is 
covariant for the gauge action on $C^*(F)$. For example, it suffices 
to set 
$$ V_z\xi:=\left\{ \begin{array}{ll} U_z\xi & \text{if } \; \xi
   \in\H_E \\ z^{-i}U_z\xi & \text{if } \; \xi\in Q_{v_i}\H_T=\H_i 
   \cong P_w\H_E.\end{array} \right. $$ 
The same argument works for graphs 
with sources: add heads as in \cite[\S1]{rs} and set
$V_z\xi:=z^{i}U_z\xi$ instead. 
\end{proof} 

We finish this preliminary section by recalling the basic definitions 
and notation about paths in a directed graph $E$.  If
$\alpha_1,\ldots,\alpha_n$ are  (not necessarily  distinct) edges such that
$r(\alpha_i)= s(\alpha_{i+1})$ for $i=1,\ldots,n-1$, then $\alpha=(\alpha_1,
\ldots,\alpha_n)$ is a \emph{path} of length $|\alpha|=n$, with source
$s(\alpha)=s(\alpha_1)$ and range  $r(\alpha)=
r(\alpha_n)$. The set of paths of length $n$ 
is denoted by $E^n$, 
$E^*:=\bigcup_{n=0}^\infty E^n$ (so that vertices in $E^0$ are identified with
paths  of length $0$), and  the set of infinite paths is denoted
$E^\infty$.  A \emph{loop} is a path of positive length whose source 
and range coincide. A loop $\alpha$ 
\emph{has an exit} if there exist an edge $e\in E^1$ and index $i$
such that $s(e)=s(\alpha_i)$ but $e\neq\alpha_i$. A graph is said to  satisfy
Condition (K) if every vertex $v\in E^0$ either lies on no loops, or there  are
two loops
$\alpha,\beta$ such that $s(\alpha)=s(\beta)=v$ and neither 
$\alpha$ nor $\beta$ is an initial subpath of the other \cite{kprr}. 

\section{Gauge-invariant ideals} 

For a row-finite graph $E$, the gauge-invariant ideals in the graph 
algebra $C^*(E)$ are in one-to-one correspondence with the 
saturated hereditary subsets of $E^0$ \cite[Theorem 4.1]{bprs}; 
indeed, if $H$ is saturated and hereditary, then the corresponding 
ideal $I_H$ is generated by $\{p_v:v\in H\}$, and if $I$ is 
a gauge-invariant ideal, then $H:=\{v\in E^0:p_v\in I\}$ is 
saturated and hereditary with $I=I_H$. When some vertices emit 
infinitely many edges, not every gauge-invariant ideal arises 
this way: in the graph 

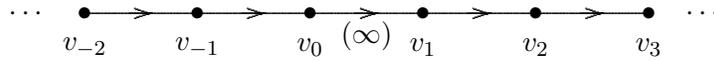
\[ \beginpicture
\setcoordinatesystem units <1.5cm,1.5cm>
\setplotarea x from -1 to 6, y from -0.5 to 0.1 

\put {$\bullet$} at 0 0
\put {$\bullet$} at 1 0
\put {$\bullet$} at 2 0
\put {$\bullet$} at 3 0 
\put {$\bullet$} at 4 0
\put {$\bullet$} at 5 0

\setlinear 
\plot 0 0 5 0 /

\arrow <0.235cm> [0.2,0.6] from 0.4 0 to 0.6 0 
\arrow <0.235cm> [0.2,0.6] from 1.4 0 to 1.6 0 
\arrow <0.235cm> [0.2,0.6] from 2.4 0 to 2.6 0 
\arrow <0.235cm> [0.2,0.6] from 3.4 0 to 3.6 0 
\arrow <0.235cm> [0.2,0.6] from 4.4 0 to 4.6 0 

\put {$\cdots$} at -0.5 0 
\put {$\cdots$} at 5.5 0 
\put {$(\infty)$} [v] at 2.5 -0.2 
\put {$v_{-2}$} [v] at 0 -0.3 
\put {$v_{-1}$} [v] at 1 -0.3 
\put {$v_{0}$} [v] at 2 -0.3 
\put {$v_{1}$} [v] at 3 -0.3 
\put {$v_{2}$} [v] at 4 -0.3 
\put {$v_{3}$} [v] at 5 -0.3 
\endpicture \] 
(where the symbol $(\infty)$ indicates infinitely many edges from 
$v_0$ to $v_1$) the projections $p_{v_i}$ associated to $H:=
\{v_i:i>0\}$ generate a gauge-invariant ideal $I_H$ with 
$H=\{v:p_v\in I_H\}$, but $H$ is not saturated in the sense 
of \cite{bprs}. So we have to adjust the notion of saturation. 

Let $E$ be a directed graph which is not necessarily row-finite. 
As usual, we write $v\geq w$ when
there is a path from $v$ to $w$, and say that a subset $H$ of $E^0$ is {\em
hereditary}  if $v\in H$ and $v\geq w$ imply $w\in H$.  Now we say that a subset
$X$ of $E^0$ is {\em saturated} if every vertex $v$ which  satisfies
$0<|s^{-1}(v)|<\infty$ and
$s(e)=v\Longrightarrow r(e)\in X$  itself belongs to $X$. (With this
definition, the set $H=\{v_i: i>0\}$ in the above example is saturated.) The
{\em saturation} 
$\Sigma(X)$ of a set $X$ is the smallest saturated set containing $X$, and
$\Sigma H(X)$ denotes the smallest saturated hereditary subset of 
$E^0$ containing $X$. 
If $v\in\Sigma(X)\setminus X$, then $0<|s^{-1}(v)|<\infty$; 
otherwise, $\Sigma(X)\setminus\{v\}$ is a smaller saturated set 
containing $X$. If $v\in\Sigma(X)$, then there is a path $\alpha$ 
with $s(\alpha)=v$ and $r(\alpha)\in X$. To see this, note 
that the elements of $\Sigma(X)$ with this property form a 
saturated set containing $X$. 

If $H$ is hereditary, so is its saturation $\Sigma(H)$. To see this, 
suppose $v\in\Sigma(H)$ and $v\geq w$, so that there is a path 
$\alpha=(\alpha_1,\ldots,\alpha_r)$ with $s(\alpha)=v$ and $r(\alpha)=w$. 
If the path enters $H$, then it stays there. So suppose $r(\alpha_i)
\not\in H$ for all $i$. Since $0<|s^{-1}(v)|<\infty$, the saturation 
property implies that $r(\alpha_1)\in\Sigma(H)$, for otherwise 
$\Sigma(H)\setminus\{r(\alpha_1)\}$ would be a smaller 
saturated set containing $H$. Since $r(\alpha_1)\in\Sigma(H)$, 
it also satisfies $0<|s^{-1}(r(\alpha_1))|<\infty$; repeating 
this argument shows that $r(\alpha_i)\in\Sigma(H)$ for all $i$, and 
in particular that $w=r(\alpha)\in\Sigma(H)$. 

\begin{remark}\label{saturation} 
For $X\subset E^0$ we can construct $\Sigma H(X)$ as the union 
of the sequence $\Sigma_n(X)$ of subsets of $E^0$ defined 
inductively as follows: 
\begin{align*}
\Sigma_0(X) & := X\cup\big\{w\in E^0:\text{ there is a path from 
   a vertex in $X$ to $w$}\big\}, \\ 
\Sigma_{n+1}(X) &:= \Sigma_n(X)\cup\big\{w\in E^0\;
:\;0<|s^{-1}(w)|<\infty \text{ and } s(e)=w \text{ imply }r(e)\in\Sigma_n(X)
\big\}. 
\end{align*} 
\end{remark} 

The next lemma provides evidence that our notion of 
saturation is the right one. 

\begin{lemma}\label{HI} 
Suppose $E$ is a directed graph and $I$ is an ideal in $C^*(E)$. Then 
$$ H_I:=\{v\in E^0:p_v\in I\} $$ 
is a saturated hereditary subset of $E^0$. 
\end{lemma} 
\begin{proof} 
Suppose first that $v\in H_I$ and $v\geq w$. Then there is a path 
$\alpha$ with $s(\alpha)=v$ and $r(\alpha)=w$, and then $s_\alpha^* 
s_\alpha=p_w$, $s_\alpha s_\alpha^*\leq p_v$. So $p_v\in I$ 
implies $s_\alpha s_\alpha^*\in I$, and $p_w=s_\alpha^*
(s_\alpha s_\alpha^*)s_\alpha\in I$. Thus $H_I$ is hereditary. 
If $v\in E^0$ has $0<|s^{-1}(v)|<\infty$ and $p_{r(e)}\in I$ 
for all $e\in E^1$ with $s(e)=v$, then $s_e=s_e p_{r(e)}\in I$ 
for all $e\in E^1$ with $s(e)=v$, and the Cuntz-Krieger 
relation (G3) at $v$ implies that $p_v\in I$; thus $H_I$ is 
saturated. 
\end{proof} 

For a hereditary subset $H$ of $E^0$, we let $I_H$ be the ideal 
of $C^*(E)$ generated by $\{p_v:v\in H\}$; since Lemma 
\ref{HI} implies that $\{v\in E^0:p_v\in I_H\}$ is 
saturated and contains $H$, we immediately have that 
$I_H=I_{\Sigma(H)}$. Since the projections generating $I_H$ 
are fixed by the gauge action it follows that the ideal 
itself is gauge-invariant. As in \cite[Lemma 4.3]{bprs}, 
we can verify that 
\begin{equation} 
I_H=\overline{\text{span}}\{s_\alpha p_v s_\beta^*:
\alpha,\beta\in E^*,\:v\in H,\:r(\alpha)=r(\beta)=v\}.  
\end{equation} 

Suppose $H$ is a saturated hereditary subset of $E^0$. When $E$ is 
row-finite, the ideal $I_H$ is Morita equivalent to the $C^*$-algebra 
$C^*(H)$ of the graph $(H,s^{-1}(H),r,s)$, and the quotient 
$C^*(E)/I_H$ is naturally isomorphic to the $C^*$-algebra of the 
graph $F:=(E^0\setminus H,r^{-1}(E^0\setminus H),r,s)$ 
\cite[Theorem 4.1]{bprs}. In general, to realise the quotient 
$C^*(E)/I_H$ as a graph algebra we have to add extra vertices to 
$F$. The problem occurs when a vertex $v$ sends infinitely many 
edges into $H$ but also finitely many into $E^0\setminus H$, 
in which case the image of the projection 
\begin{equation} 
p_{v,H}:=\sum_{s(e)=v,\;r(e)\not\in H}s_e s_e^* 
\end{equation} 
will be strictly smaller in $C^*(E)/I_H$ than the image of $p_v$. 
To get round this, we add a new 
sink $\beta(v)$ to $F^0$ and extra edges $\beta(e)$ with 
$r(\beta(e))=\beta(v)$ for each edge $e$ with $r(e)=v$. 

Formally, we define $\Hbad$ to be the set of such 
vertices; thus 
$$ \Hbad:=\{v\in E^0\setminus H:|s^{-1}(v)|=\infty 
   \text{ and } 0<|s^{-1}(v)\cap r^{-1}(E^0\setminus H)|<\infty\}. $$ 
We then define a graph $E/H$ by 
\begin{align*}
(E/H)^0 & := (E^0\setminus H)\cup\{\beta(v):v\in \Hbad \}, \\ 
(E/H)^1 & := r^{-1}(E^0\setminus H)\cup\{\beta(e):e\in E^1,r(e)\in 
              \Hbad\},  
\end{align*} 
with $r,s$ extended by $s(\beta(e))=s(e)$ and $r(\beta(e))=\beta(r(e))$. 

\begin{example}\label{quotientgraph} 
In the following graph,  

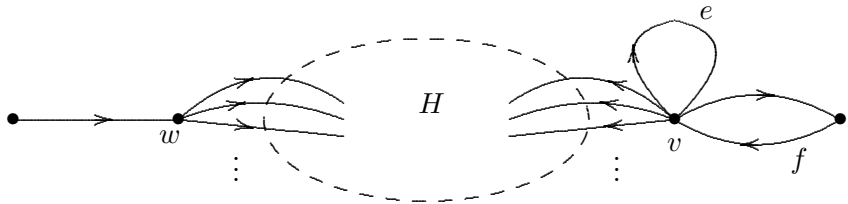
\[ \beginpicture
\setcoordinatesystem units <1.1cm,1.1cm>
\setplotarea x from -6 to 1.7, y from -1 to 0.5
\put {$\bullet$} at -6 0
\put {$\bullet$} at -4 0
\put {$\bullet$} at 2 0
\put {$\bullet$} at 4 0
\setlinear 
\plot -6 0 -4 0 /
\setquadratic
\plot -4 0 -3 0.5 -2 0.2 /
\plot -4 0 -3 0.2 -2 0 /
\plot -4 0 -3 -0.1 -2 -0.2 /
\plot 2 0 1 0.5 0 0.2 /
\plot 2 0 1 0.2 0 0 /
\plot 2 0 1 -0.1 0 -0.2 /
\plot 2 0 3 0.3 4 0 /
\plot 2 0 3 -0.3 4 0 /
\plot 2 0 1.5 0.8 2 1.2 2.5 0.8 2 0 /
\put {$H$} [l] at -1.1 0.2
\put {$v$} [v] at 2 -0.3 
\put {$w$} [v] at -4.1 -0.2 
\put {$e$} [l] at 2.3 1.3
\put {$f$} [v] at 3.5 -0.5
\arrow <0.235cm> [0.2,0.6] from -5 0 to -4.8 0 
\arrow <0.235cm> [0.2,0.6] from -3.2 0.47 to -3.1 0.5 
\arrow <0.235cm> [0.2,0.6] from -3.3 0.19 to -3.2 0.21
\arrow <0.235cm> [0.2,0.6] from -3.2 -0.08 to -3.1 -0.1
\arrow <0.235cm> [0.2,0.6] from 1.25 0.455 to 1.2 0.47
\arrow <0.235cm> [0.2,0.6] from 1.2 0.2 to 1.1 0.22
\arrow <0.235cm> [0.2,0.6] from 1.3 -0.07 to 1.15 -0.08
\arrow <0.235cm> [0.2,0.6] from 3 0.3 to 3.2 0.3 
\arrow <0.235cm> [0.2,0.6] from 3 -0.3 to 2.8 -0.3
\arrow <0.235cm> [0.2,0.6] from 1.5 0.7 to 1.501 0.9  
\setdashes 
\ellipticalarc axes ratio 2:1 360 degrees from 0.9 .25 center at -1 0
\put {$\vdots$} at -3.3 -0.5
\put {$\vdots$} at 1.3 -0.5 
\endpicture \] 
we have $\Hbad=\{v\}$, and the graph $E/H$ looks like 
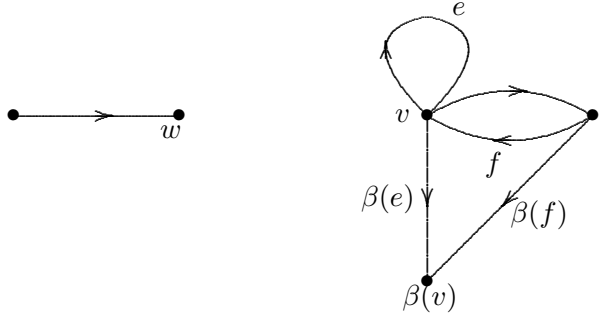
\[ \beginpicture
\setcoordinatesystem units <1.1cm,1.1cm>
\setplotarea x from -6 to 10, y from -2 to 1
\put {$\bullet$} at -4 0
\put {$\bullet$} at -2 0
\put {$\bullet$} at 1 0
\put {$\bullet$} at 3 0
\put {$\bullet$} at 1 -2
\setlinear 
\plot -4 0 -2 0 / 
\plot 1 0 1 -2 /
\plot 3 0 1 -2 /
\setquadratic
\plot 1 0 2 0.3 3 0 /
\plot 1 0 2 -0.3 3 0 /
\plot 1 0 0.5 0.8 1 1.2 1.5 0.8 1 0 /
\put {$v$} [v] at 0.7 0 
\put {$w$} [v] at -2.1 -0.2 
\put {$e$} [l] at 1.3 1.3
\put {$f$} [v] at 1.8 -0.6
\put {$\beta(v)$} [l] at 0.7 -2.2 
\put {$\beta(e)$} [l] at 0.2 -1 
\put {$\beta(f)$} [r] at 2.7 -1.2
\arrow <0.235cm> [0.2,0.6] from -3 0 to -2.8 0 
\arrow <0.235cm> [0.2,0.6] from 1 -1 to 1 -1.1
\arrow <0.235cm> [0.2,0.6] from 2 -1 to 1.9 -1.1
\arrow <0.235cm> [0.2,0.6] from 2 0.3 to 2.2 0.3 
\arrow <0.235cm> [0.2,0.6] from 2 -0.3 to 1.8 -0.3
\arrow <0.235cm> [0.2,0.6] from 0.5 0.7 to 0.501 0.9  
\endpicture \] 
\end{example} 

\begin{prop}\label{quotientalgebra} 
Let $H$ be a hereditary subset of a directed graph $E$. Then the ideal 
$I_H$ defined in \textnormal{(1)} is Morita equivalent to the $C^*$-algebra of
the  graph $(H,s^{-1}(H),r,s)$. Let $\pi:C^*(E)\rightarrow C^*(E)/I_H$ be the 
quotient map, let $\{s_e,p_v\}$ be the canonical Cuntz-Krieger $E$-family, 
and write $S_e=\pi(s_e)$, $P_v=\pi(p_v)$, $P_{v,H}=\pi(p_{v,H})$,  
where $p_{v,H}$ are the projections defined in \textnormal{(2)}. 
If $H$ is also saturated, then
\begin{align} 
\mbox{\hskip1in}Q_v & := P_v &&
\text{if $v\in(E/H)^0\setminus\Hbad$}\mbox{\hskip1in} \notag\\  
Q_v & := P_{v, H}&&\text{if $v\in \Hbad$} \notag\\  
Q_{\beta(v)} & := P_v-P_{v,H}&&\text{if $v\in \Hbad$}\\
T_e & := S_e&&\text{if $r(e)\in(E^0\setminus H)
\setminus\Hbad$}\notag\\ 
T_e & := S_e P_{r(e),H}&&\text{if $r(e)\in(E^0\setminus H)\cap\Hbad$}\notag\\
T_{\beta(e)} & := S_e (P_{r(e)}-P_{r(e),H})&&\text{if $r(e)\in(E^0\setminus
H)\cap\Hbad$}\notag
\end{align}
is a Cuntz-Krieger $(E/H)$-family in $C^*(E)/I_H$, and the homomorphism 
$\pi_{T,Q}$ is an isomorphism of $C^*(E/H)$ onto $C^*(E)/I_H$. 
\end{prop} 
\begin{proof} 
The argument of \cite[Theorem 4.1(c)]{bprs} shows that there is a natural 
isomorphism of $C^*(H)$ onto the corner of $I_H$ determined by the projection 
$p_H:=\sum_{v\in H}p_v\in M(I_H)$, and that this projection is full. 

It is tedious but straightforward to verify that $\{T_e,Q_v\}$ 
is a Cuntz-Krieger $(E/H)$-family, and hence there is a homomorphism 
$\pi_{T,Q}:C^*(E/H)\rightarrow C^*(E)/I_H$ which carries the generating 
family $\{t_e,q_v\}$ of $C^*(E/H)$ into $\{T_e,Q_v\}$. To see that 
$\pi_{T,Q}$ is surjective, note that we can recover $\{S_e,P_v\}$ 
from $\{T_e,Q_v\}$: 
\begin{equation} 
\begin{split} 
P_v & = \begin{cases} Q_v & \text{if $v\not\in H\cup
          \Hbad$} \\ Q_v+Q_{\beta(v)} & \text{if $v\in \Hbad$} \\ 
          0 & \text{if $v\in H$} \end{cases} \\ 
S_e & = \begin{cases} T_e & \text{if $r(e)\not\in 
          H\cup\Hbad$} \\ T_e+T_{\beta(e)} & \text{if $r(e)\in 
          \Hbad$} \\ 0 & \text{if $r(e)\in H$} \end{cases}  
\end{split}  
\end{equation} 
The formulas (4) also show how to construct a Cuntz-Krieger $E$-family 
$\{S_e,P_v\}$ from a Cuntz-Krieger $(E/H)$-family $\{T_e,Q_v\}$ in such 
a way that the formulas (3) recover $\{T_e,Q_v\}$. Thus there are 
Cuntz-Krieger $E$-families $\{S_e,P_v\}$ with $P_v=0$ for $v\in H$ 
such that the projections $Q_v,Q_{\beta(v)}$ defined in (3) are all 
non-zero, and in particular this must be true of those defined by 
the universal family $\{\pi(s_e),\pi(p_v)\}$. It therefore follows 
from gauge-invariant uniqueness (Theorem~\ref{uniqueness})  
that the homomorphism $\pi_{T,Q}$ is injective. 
\end{proof} 

The sinks in a directed graph give rise to ideals: indeed, if $B\subset E^0$ 
consists of sinks, then the projections $p_v$ associated to the sinks 
generate a family $\{I_v:v\in B\}$ of mutually orthogonal 
ideals. Since the vertices $\{\beta(v):v\in \Hbad\}$ 
are sinks in $E/H$, they give rise to ideals in $C^*(E/H)\cong 
C^*(E)/I_H$, and hence to ideals in $C^*(E)$. More formally, 
if $H$ is saturated and hereditary,  
then for $B\subset \Hbad$ we let $J_{H,B}$ denote the 
ideal of $C^*(E)$ generated by the projections $\{p_v:v\in H\}\cup
\{p_v-p_{v,H}:v\in B\}$. The usual arguments show that 
$$ J_{H,B}=\overline{\text{span}}\big\{s_\alpha p_v s_\beta^*,\;
   s_\mu(p_w-p_{w,H})s_\nu^*:v\in H,\ \alpha,\beta\in 
r^{-1}(v),\ w \in B,\ 
   \mu,\nu\in r^{-1}(w)\big\} $$ 
and that $J_{H,B}$ is gauge-invariant. Note also that 
$I_H=J_{H,\emptyset}$, and that $I_H\subset J_{H,B}$ for all $B$. 

To identify the quotient $C^*(E)/J_{H,B}$, note that the set 
$\beta(B)$ is saturated and hereditary in $(E/H)^0$. Since 
the quotient map $C^*(E)\rightarrow C^*(E/H)=C^*(t_e,q_v)$ takes 
$p_v-p_{v,H}$ into $q_{\beta(v)}$, it maps $J_{H,B}$ onto 
the ideal $I_{\beta(B)}$ of $C^*(E/H)$ generated by 
$\{q_{\beta(v)}:v\in B\}$. Since $\beta(B)^{\text{fin}}_\infty=
\emptyset$ in $E/H$, the second quotient $(E/H)/\beta(B)$ is just 
$(E/H)\setminus\beta(B)=\left( (E/H)^0\setminus\beta(B),
r^{-1}((E/H)^0\setminus\beta(B)),r,s\right)$. Thus we have: 

\begin{cor}\label{secondquotient} 
If $H$ is a saturated hereditary subset of $E^0$ and 
$B\subset \Hbad$, then $C^*(E)/J_{H,B}$ is naturally 
isomorphic to $C^*((E/H)\setminus\beta(B))$. 
\end{cor} 

The following theorem gives a complete list of the gauge-invariant 
ideals of $C^*(E)$ for an arbitrary infinite graph $E$. 

\begin{thm}\label{gaugeinvariant} 
Let $E$ be a directed graph. Then the ideals 
$$ \big\{ J_{H,B}:H \text{ is saturated and hereditary, } 
   B\subset\Hbad\big\} $$ 
are distinct gauge-invariant ideals in $C^*(E)$, and every 
gauge-invariant ideal is of this form. Indeed, if $I$ is 
a gauge-invariant ideal in $C^*(E)=C^*(s_e,p_v)$, $H:=\{
v\in E^0:p_v\in I\}$, and $B:=\{v\in\Hbad:p_v-p_{v,H}\in I\}$, 
then $I=J_{H,B}$. 
\end{thm} 

We begin by showing that we can recover $H$ and $B$ from $J_{H,B}$; 
this immediately implies that the ideals are distinct. 

\begin{lemma}\label{HB} 
Let $E$ be a directed graph. Suppose that $H$ is a saturated 
hereditary subset of $E^0$ and $B\subset \Hbad$. Then 
$H=H_{J_{H,B}}$. If we use the isomorphism of Proposition 
\ref{quotientalgebra} to view $J_{H,B}/I_H$ as an ideal in 
$C^*(E/H)$, then $\beta(B)=H_{J_{H,B}/I_H}$. 
\end{lemma} 
\begin{proof} 
We trivially have $H\subset H_{J_{H,B}}$. Suppose $v\not\in H$. Then 
the image of $p_v$ under the isomorphism of $C^*(E)/J_{H,B}$ onto 
$C^*((E/H)\setminus\beta(B))=C^*(t_e,q_v)$ dominates the 
projection $q_v$ associated to the vertex $v\in((E/H)\setminus
\beta(B))^0$, and hence is nonzero; thus $p_v\not\in J_{H,B}$. 
This gives the first assertion. For the second, note that the 
image of $J_{H,B}$ under the quotient map $C^*(E)\rightarrow 
C^*(E/H)=C^*(u_f,r_w)$ is the ideal $I_{\beta(B)}$, and 
$\beta(B)=\{w\in(E/H)^0:r_w\in I_{\beta(B)}\}$ by 
the first assertion. 
\end{proof} 

\begin{proof}[Proof of Theorem \ref{gaugeinvariant}] 
Lemma \ref{HB} implies that the ideals are distinct. Given $I$ and 
$H,B$ as in the theorem, we note that $J_{H,B}\subset I$, and 
consider the image of $I/J_{H,B}$ in $C^*((E/H)\setminus\beta(B))
=C^*(t_e,q_v)$. We shall show by contradiction that there is no vertex 
$w$ of $(E/H)\setminus\beta(B)$ such that the corresponding projection 
$q_w$ lies in $I/J_{H,B}$. If $w\in E^0\setminus H$ and $w\not\in 
\Hbad$, then $q_w\in I/J_{H,B}$ implies $p_w\in I$, which 
contradicts $w\not\in H$. If $w\in \Hbad$, then $q_w\in 
I/J_{H,B}$ implies $p_{w,H}\in I$; now we can choose $e\in E^1$ 
such that $s(e)=w$ and $r(e)\not\in H$, and then $p_{w,H}\in I$ 
implies $p_{r(e)}=s_e s_e^*\in I$, which is incompatible with 
$r(e)\not\in H$. If $w=\beta(v)$ for some $v\in \Hbad$, then 
$q_w\in I/J_{H,B}$ implies $p_v-p_{v,H}\in I$, and $w=\beta(v)
\in\beta(B)$, which is impossible because $w$ is a vertex of 
$(E/H)\setminus\beta(B)$. Thus for all $w\in((E/H)\setminus\beta(B))^0$, 
$q_w$ has non-zero image in $C^*((E/H)\setminus\beta(B))/
(I/J_{H,B})$. Now gauge-invariant uniqueness
(Theorem~\ref{uniqueness}) implies  that the quotient map of
$C^*((E/H)\setminus\beta(B))$ onto 
$C^*((E/H)\setminus\beta(B))/(I/J_{H,B})$ is injective, which says that 
$I/J_{H,B}=0$ and $I=J_{H,B}$. 
\end{proof}

\begin{cor}\label{conditionK} 
Suppose $E$ is a directed graph satisfying Condition~\textnormal{(K)}. Then
every ideal of $C^*(E)$ is gauge-invariant, and hence
Theorem~\ref{gaugeinvariant} gives a complete description of the
ideals  of
$C^*(E)$. 
\end{cor}

\begin{proof}
Suppose $I$ is an ideal in $C^*(E)$. Let $H:=\{v\in E^0:p_v\in I\}$, which
is saturated and hereditary by Lemma~\ref{HI}, and let
$B:=\{v\in\Hbad:p_v-p_{v,H}\in I\}$. Note that $J_{H,B}\subset I$. Let
$(E/H)\setminus\beta(B)$ denote the graph appearing in
Corollary~\ref{secondquotient}. As in the first paragraph of the proof of
\cite[Theorem~4.1]{bprs}, both quotients $C^*(E)/I$ and $C^*(E)/J_{H,B}$ are
generated by Cuntz-Krieger $((E/H)\setminus\beta(B))$-families in which all the
projections associated to vertices are nonzero.

We claim that all loops in $(E/H)\setminus\beta(B)$ have exits. Suppose
$\alpha$ is a loop in $(E/H)\setminus\beta(B)$. Since all the new vertices
added to $E^0\setminus H$ to form $(E/H)\setminus\beta(B)$ are sinks, the loop
$\alpha$ must come from a loop $\tilde\alpha$ in $E$. Because
$E$ satisfies (K), each vertex in $\tilde\alpha$ must lie on another loop.
Since this loop cannot enter the hereditary set $H$, there must be an exit from
$\tilde\alpha$ which lies in $r^{-1}(E^0\setminus H)$, and hence gives an exit
from $\alpha$ in $(E/H)\setminus\beta(B)$. This justifies the claim.

Now two applications of the Cuntz-Krieger uniqueness theorem
\cite[Theorem~1.5]{rs} show that both quotients $C^*(E)/I$ and $C^*(E)/J_{H,B}$
are canonically isomorphic to $C^*((E/H)\setminus\beta(B))$. Thus the
quotient map of $C^*(E)/J_{H,B}$ onto $C^*(E)/I$ is an isomorphism, and
$I=J_{H,B}$. The corollary now follows from Theorem~\ref{gaugeinvariant}.
\end{proof} 

We need the following proposition  in the proof of Lemma~\ref{maximaltail}
below, and in the analysis of the hull-kernel  topology on the primitive
ideal space in our sequel. 

\begin{prop}\label{intersection} 
Suppose $E$ is a directed graph, $\{H_i:i\in\Lambda\}$ is 
a family of saturated hereditary subsets of $E^0$, and $B_i\subset
\Hibad$ for $i\in\Lambda$. Let $H=\bigcap_{i\in\Lambda}H_i$ and 
$B=\left(\bigcap_{i\in\Lambda}H_i\cup B_i\right)\cap\Hbad$. Then 
$$ \textstyle{\bigcap_{i\in\Lambda}J_{H_i,B_i}=J_{H,B}.} $$ 
\end{prop} 
\begin{proof} 
Since the intersection of gauge-invariant ideals is gauge-invariant, 
Theorem \ref{gaugeinvariant} says that $\bigcap_{i\in\Lambda}J_{H_i,B_i}
=J_{K,C}$ for 
$$ \textstyle{K=\big\{v:p_v\in\bigcap_{i\in\Lambda}J_{H_i,B_i}\big\} \; \text{
and } 
   \; C=\big\{w\in\Kbad:p_w-p_{w,K}\in\bigcap_{i\in\Lambda}
   J_{H_i,B_i}\big\}.} $$ 
By two applications of Lemma \ref{HB}, we have 
$$ \textstyle{K=H_{J_{K,C}}=\bigcap_{i\in\Lambda}H_{J_{H_i,B_i}}=
   \bigcap_{i\in\Lambda}H_i=H.} $$ 
It remains to identify $C$ with $B$. Let $w\in\Kbad$; we want to show that 
$w\in\bigcap_{i\in\Lambda}H_i\cup B_i$ if and only if $w\in C$. 

Suppose $w\in C$ and $i\in\Lambda$ is fixed. Then $p_w-p_{w,K}\in J_{K,C}
\subset J_{H_i,B_i}$. For each of the finitely many $e$ with $s(e)=w$ and 
$r(e)\not\in K$, $r(e)\in H_i$ implies $p_{r(e)}\in J_{H_i,B_i}$ and 
$s_e s_e^*\in J_{H_i,B_i}$. If $r(e)\in H_i$ for all such $e$, then 
$p_{w,K}\in J_{H_i,B_i}$, $p_w=(p_w-p_{w,K})+p_{w,K}\in J_{H_i,B_i}$, 
and $w\in H_i$. If $r(e)\not\in H_i$ for some such $e$, then 
$$ p_w-p_{w,H_i}=p_w-p_{w,K}+\sum_{s(e)=w,\;r(e)\in H_i\setminus K}
   s_e s_e^*\in J_{H_i,B_i}, $$ 
and $w\in B_i$. Either way, $w\in H_i\cup B_i$. 

For the converse, suppose $w\in H_i\cup B_i$ for all $i$, and fix $i$; 
we want to show $p_w-p_{w,K}\in J_{H_i,B_i}$. If $w\in H_i$, then 
$p_w\in J_{H_i,B_i}$, so this is trivially true. If $w\in B_i$, then 
$p_w-p_{w,H_i}\in J_{H_i,B_i}$, and 
$$ p_w-p_{w,K}=p_w-p_{w,H_i}+\sum_{s(e)=w,\;r(e)\in H_i\setminus K}
   s_e s_e^*\in J_{H_i,B_i}, $$ 
as required. 
\end{proof} 

\begin{cor}\label{inclusion} 
If $E$ is a directed graph, $H_1$ and $H_2$ are saturated hereditary 
subsets of $E^0$, and $B_i\subset \Hibad$ for $i=1,2$, then 
$J_{H_1,B_1}\subset J_{H_2,B_2}$ if and only if 
$H_1\subset H_2$ and $B_1\subset H_2\cup B_2$.  
\end{cor} 

\section{Gauge-invariant primitive ideals} 

The primitive ideal space of the $C^*$-algebras of row-finite 
graphs satisfying Condition (K) were described in \cite[\S6]{bprs}. 
In particular, \cite[Corollary 6.5]{bprs} gives a bijection between 
the primitive ideals and certain subsets of the vertex set, called 
maximal tails. The concept of a maximal tail also plays a crucial role 
in our analysis of primitive gauge-invariant ideals in $C^*(E)$. However,
we need to adjust the definition to accommodate non-row-finite graphs. 

\begin{lemma}\label{maximaltail}
Suppose $I$ is an ideal in $C^*(E)$. Then $M:=E^0\setminus H_I$ 
satisfies: 
\begin{itemize} 
\item[(a)] if $v\in E^0$, $w\in M$, and $v\geq w$ in $E$, then $v\in M$, and 
\item[(b)] if $v\in M$ and $0<|s^{-1}(v)|<\infty$, then there exists 
$e\in E^1$ with $s(e)=v$ and $r(e)\in M$. 
\end{itemize} 
If $I$ is a primitive ideal, then in addition: 
\begin{itemize} 
\item[(c)] for every $v,w\in M$ there exists $y\in M$ such that 
$v\geq y$ and $w\geq y$. 
\end{itemize} 
\end{lemma} 
\begin{proof} 
Conditions (a) and (b) say that $H_I$ is hereditary and saturated, 
so the first part follows from Lemma \ref{HI}. 

Now suppose $I$ is primitive, $v,w\in M$, and there is no $y$ 
as in (c). The sets $H_v=\{x\in E^0:v\geq x\}$ and $H_w=\{x\in 
E^0:w\geq x\}$ are hereditary and satisfy $H_v\cap H_w\subset H_I$.  
The set $H_I\cup(E^0\setminus H_w)$ is saturated. 
Indeed, let $x\in E^0$ be such that $0<|s^{-1}(x)|<\infty$ and 
$r(e)\in H_I\cup(E^0\setminus H_w)$ for each edge $e$ with $s(e)=x$. 
Then either $x\in H_I$ or there is at least one $e$ with $s(e)=x$ and 
$r(e)\not\in H_w$, in which case $x\in E^0\setminus H_w$ because 
$H_w$ is hereditary. Thus $\Sigma H(v)=\Sigma(H_v)\subset H_I\cup 
(E^0\setminus H_w)$. The same argument shows $\Sigma H(w)\subset 
H_I\cup(E^0\setminus H_v)$, and the hypothesis $H_v\cap H_w\subset H_I$ forces
$\Sigma H(v)\cap\Sigma H(w)
\subset H_I$. It follows from Lemma \ref{intersection} that 
$$ I_{\Sigma H(v)}\cap I_{\Sigma H(w)}= 
   I_{\Sigma H(v)\,\cap\,\Sigma H(w)}\subset I, $$ 
which is impossible because neither $I_{\Sigma H(v)}$ nor 
$I_{\Sigma H(w)}$ is contained in $I$. Thus there must exist
$y\in M$ such that $v\geq y$ and $w\geq y$. 
\end{proof} 

We now define a {\em maximal tail} in $E$ to be a nonempty subset 
$M$ of $E^0$ satisfying conditions (a), (b) and (c) of Lemma 
\ref{maximaltail}. If $M$ is a maximal tail in $E$ then we say that 
{\em every loop in $M$ has an exit} if every loop with vertices 
in $M$ has an exit $e\in E^1$ with $r(e)\in M$. 

\begin{prop}\label{criterionprimitive} 
Let $E$ be a directed graph. Then $C^*(E)$ is primitive if and 
only if every loop in $E$ has an exit and for every $v,w\in E^0$ there 
exists $y\in E^0$ such that $v\geq y$ and $w\geq y$. 
\end{prop} 
\begin{proof} 
Suppose that $C^*(E)$ is primitive. Then $E^0$ is a maximal tail 
by Lemma \ref{maximaltail}. If there is a loop $L$ in $E$ without 
exits and $L^0$ is the set of vertices on
$L$, then the ideal $I_{L^0}$ is Morita equivalent to
$C(\T)$  (cf. \cite[\S2]{kpr} and formula (1) above), contradicting 
primitivity of $C^*(E)$. 

Conversely, suppose the two conditions of the proposition are satisfied. 
It suffices to show that $\{0\}$ is a prime ideal of $C^*(E)$. 
Indeed, let $I_1$, $I_2$ be two non-zero ideals of $C^*(E)$. Since 
all loops in $E$ have exits, \cite[Theorem 2]{flr} implies that 
there exist $v,w\in E^0$ such that $p_v\in I_1$ and 
$p_w\in I_2$. If $y\in E^0$ satisfies $v\geq y$ and $w\geq y$, 
then we have $p_y\in I_1\cap I_2$, so $I_1\cap I_2\not=\{0\}$, and $\{0\}$ is
prime.  
\end{proof} 

\begin{remark}\label{primitivegraph}
The conditions of Proposition  
\ref{criterionprimitive} are equivalent to:
\begin{enumerate} 
\item for every $v,w\in E^0$ we have 
$\Sigma H(v)\cap\Sigma H(w)\neq\emptyset$, and 
\item every loop in $E$ has an exit. 
\end{enumerate} 
The above two conditions should be compared with the 
criterion of simplicity for graph algebras \cite[Theorem 12]{sz1}: 
\begin{enumerate} 
\item for every $v\in E^0$ we have $\Sigma H(v)=E^0$, and 
\item every loop in $E$ has an exit. 
\end{enumerate} 
\end{remark} 

\begin{cor}\label{primitivequotient} 
Let $E$ be a directed graph and let $H$ be a saturated 
hereditary subset of $E^0$. Then $C^*(E\setminus H)$ is primitive 
if and only if $M:=E^0\setminus H$ is a maximal tail such that all 
loops in $M$ have exits. 
\end{cor}

By Theorem \ref{gaugeinvariant}, we can determine 
all the gauge-invariant primitive ideals of $C^*(E)$ by deciding which of the
ideals $J_{H,B}$ are primitive. To this  end, we use Corollary
\ref{primitivequotient} to see which  quotient algebras $C^*(E)/J_{H,B}\cong
C^*((E/H)\setminus
\beta(B))$ are primitive. If $\Hbad\setminus B$ contains distinct 
vertices $v,w$, there are at least two sinks $\beta(v),\beta(w)$ 
in $(E/H)\setminus\beta(B)$ and $C^*((E/H)\setminus\beta(B))$ 
cannot be primitive by Proposition \ref{criterionprimitive}. 
So we only need to consider two possibilities: 
$B=\Hbad$ and $B=\Hbad\setminus\{v\}$ for some $v\in\Hbad$. 

\begin{lemma}\label{none} 
Let $E$ be a directed graph and let $H$ be a saturated and 
hereditary subset of $E^0$. Then $J_{H,\Hbad}$ is primitive 
if and only if $M:=E^0\setminus H$ is a maximal tail such that  
all loops in $M$ have exits. 
\end{lemma} 
\begin{proof} 
Since $C^*(E)/J_{H,\Hbad}$ is isomorphic to $C^*((E/H)\setminus\beta(\Hbad))=
C^*(E\setminus H)$, this follows from Corollary \ref{primitivequotient}. 
\end{proof} 

For any non-empty subset $X$ of $E^0$ we denote by $\Omega(X)$ 
the collection of vertices $w\in E^0\setminus X$ 
such that there is no path from $w$ to any vertex in $X$. That is, 
$$ \Omega(X):=\left\{w\in E^0\setminus X:
   w\not\geq v \mbox{ for all $v\in X$}\right\}. $$ 
Note that we have $\Omega(M)=E^0\setminus M$ 
for any maximal tail $M$. The following lemma shows that it is important 
to look at the sets $\Omega(v)$ corresponding to certain
vertices $v$. 

\begin{lemma}\label{one} 
Let $E$ be a directed graph, let $H$ be a saturated 
hereditary subset of $E^0$, and let $v\in\Hbad$. Then 
$J_{H,\Hbad\setminus\{v\}}$ is primitive if and only if $H=\Omega(v)$. 
\end{lemma} 
\begin{proof} 
The ideal $J_{H,\Hbad\setminus\{v\}}$ is primitive if and only if 
the corresponding quotient algebra $C^*((E/H)\setminus\beta
(\Hbad\setminus\{v\}))$ is primitive. Since the 
graph $(E/H)\setminus\beta(\Hbad\setminus\{v\})$ 
contains a sink $\beta(v)$, Corollary \ref{primitivequotient} implies 
that $J_{H,\Hbad\setminus\{v\}}$ is primitive if and only if 
for any vertex $w\in E^0\setminus H$ there exists a path from $w$ to $v$. 
This, however, is equivalent to $H=\Omega(v)$. 
\end{proof} 

We call a vertex $v\in E^0$ with the property described 
in Lemma~\ref{one} a {\em breaking vertex}, and write $v\in BV(E)$. More 
formally,  we define
$$ BV(E):=\left\{v\in E^0:|s^{-1}(v)|=\infty \text{ and } 
   0<|s^{-1}(v)\setminus\Omega(v)|<\infty\right\}. $$ 
Note that if a vertex $v$ emits infinitely many edges then 
$\Omega(v)$ is automatically saturated and hereditary, and hence 
$v$ is a breaking vertex if and only if $v\in\Obad$. If the graph $E$ is
row-finite, there are no breaking vertices. 

\begin{thm}\label{primitivegaugeinvariant} 
Let $E$ be a directed graph. Then the gauge-invariant 
primitive ideals in $C^*(E)$ are the ideals $J_{\Omega(M),
\Omega(M)_\infty^{\text{{\rm fin}}}}$ 
associated to the maximal tails $M$ in which all loops have exits, 
and the ideals $J_{\Omega(v),\Obad\setminus\{v\}}$ associated 
to breaking vertices $v\in BV(E)$. These ideals are distinct. 
\end{thm} 
\begin{proof} 
By Theorem \ref{gaugeinvariant}, all gauge-invariant ideals in 
$C^*(E)$ have the form $J_{H,B}$, with $H$ a saturated hereditary 
subset of $E^0$ and $B$ a subset of $\Hbad$, and these 
ideals are distinct. So we only need to decide which of these 
ideals are primitive. If $\Hbad\setminus B$ has two or more 
vertices then $J_{H,B}$ is not primitive, since $C^*(E)/J_{H,B}
\cong C^*((E/H)\setminus\beta(B))$ and the graph $(E/H)\setminus
\beta(B)$ contains at least two sinks, contradicting the conditions 
of Proposition \ref{criterionprimitive}. Thus we may assume that 
either $B=\Hbad$ or $B=\Hbad\setminus\{v\}$ for some $v\in\Hbad$. 
If $B=\Hbad$, the ideal $J_{H,B}$ is primitive if and only 
if $H=\Omega(M)$ for some maximal tail $M$ in which all loops have 
exits by Lemma \ref{none}. If $B=\Hbad\setminus\{v\}$, the ideal $J_{H,B}$ 
is primitive if and only if $H=\Omega(v)$ for some breaking vertex 
$v\in BV(E)$ by Lemma \ref{one}.  
\end{proof}

The following corollary now follows from Corollary~\ref{conditionK} and
Theorem~\ref{primitivegaugeinvariant}. 

\begin{cor}\label{conditionKprimitive} 
If $E$ is a directed graph satisfying Condition \textnormal{(K)} then Theorem 
\ref{primitivegaugeinvariant} gives a complete description 
of the primitive ideals of $C^*(E)$. 
\end{cor} 

\section{Examples} 

We illustrate our results with four examples which are not covered by the
existing literature. 

\begin{example}\label{example1} 
The following directed graph satisfies Condition (K) but is not 
row-finite. 

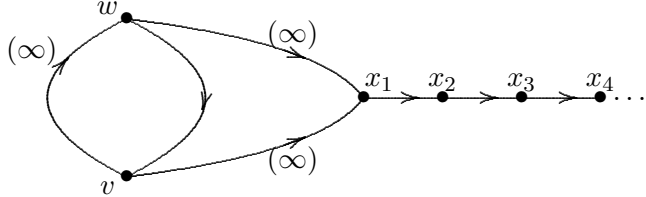
\[ \beginpicture
\setcoordinatesystem units <2.1cm,2.1cm>
\setplotarea x from -4 to 1, y from -0.3 to 0.7

\put {$\bullet$} at -3 0
\put {$\bullet$} at -3 1
\put {$\bullet$} at -1.5 0.5
\put {$\bullet$} at -1 0.5
\put {$\bullet$} at -0.5 0.5
\put {$\bullet$} at 0 0.5

\put {$\ldots$} at 0.2 0.5
\put {$(\infty)$} at -3.6 0.8
\put {$(\infty)$} at -1.95 0.9
\put {$(\infty)$} at -1.95 0.1

\setlinear 
\plot -1.5 0.5 0 0.5 / 

\arrow <0.235cm> [0.2,0.6] from -1.2 0.5 to -1.15 0.5 
\arrow <0.235cm> [0.2,0.6] from -0.7 0.5 to -0.65 0.5 
\arrow <0.235cm> [0.2,0.6] from -0.2 0.5 to -0.15 0.5 

\setquadratic 
\plot -3 1 -2 0.8 -1.5 0.5 /
\plot -3 0 -2 0.2 -1.5 0.5 /
\plot -3 0 -3.5 0.5 -3 1 /
\plot -3 1 -2.5 0.5 -3 0 /

\arrow <0.235cm> [0.2,0.6] from -2 0.8 to -1.87 0.75 
\arrow <0.235cm> [0.2,0.6] from -2 0.2 to -1.87 0.25

\arrow <0.235cm> [0.2,0.6] from -3.41 0.71 to -3.37 0.752
\arrow <0.236cm> [0.2,0.6] from -2.515 0.445 to -2.52 0.425

\put {$v$} [v] at -3.12 -0.07 
\put {$w$} [v] at -3.12  1.07 
\put {$x_1$} [v] at -1.4  0.6 
\put {$x_2$} [v] at -1 0.6 
\put {$x_3$} [v] at -0.5 0.6 
\put {$x_4$} [v] at 0 0.6 
 
\endpicture \] 
There are only two maximal 
tails, $M_1=E^0$ and $M_2=\{v,w\}$, and in both every loop has 
an exit. The  
primitive ideal corresponding to $M_1$ is $\{0\}$, and hence $C^*(E)$ is
primitive. The primitive ideal corresponding to
$M_2$ is $J_{X,\{w\}}$, where $X=\{x_1,x_2,x_3,\ldots\}$. The only breaking
vertex is $w$, and the corresponding primitive ideal is $I_X$. 

Since $I_X$ is infinite-dimensional and Morita
equivalent  to $C^*(X)$ \cite[\S2]{kpr}, and $C^*(X)$  is isomorphic to
the
$C^*$-algebra $\K$ of compact operators on a separable infinite-dimensional
Hilbert space, we have $I_X\cong \K$ also. By Proposition
\ref{quotientalgebra},  the quotient graph $E/X$ contains one sink $\beta(w)$,
and $\beta(w)$  is the range of infinitely many paths. Thus
$J_{X,\{w\}}/I_X\cong\K$.  The ideal  $J_{X,\{w\}}$ is an
extension of $\K$ by $\K$, and is the unique essential extension by \cite[Lemma
1.1]{sz4}. Another application of  Proposition \ref{quotientalgebra} shows that
$C^*(E)/J_{X,\{w\}}\cong  C^*(E\setminus X)\cong M_2(\C)\otimes\O_\infty$. 
\end{example} 

\begin{example}\label{example2} 
Let $E$ be the following graph: 

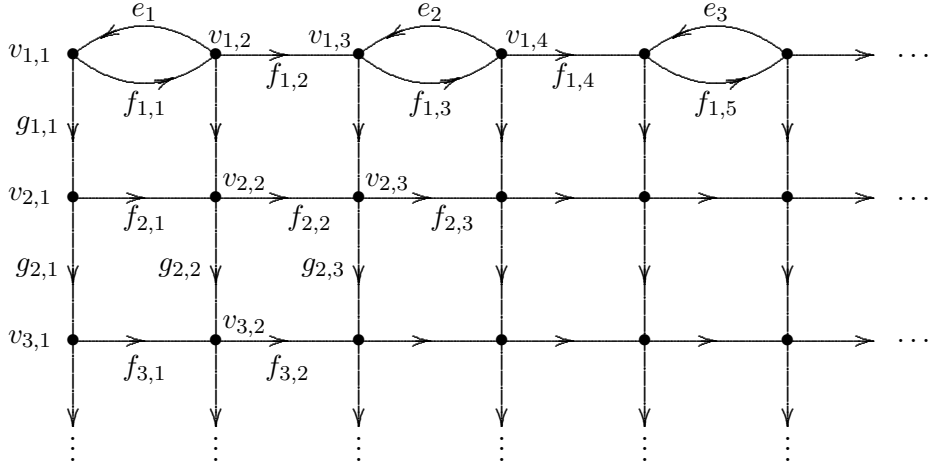
\[ \beginpicture
\setcoordinatesystem units <1.9cm,1.9cm>
\setplotarea x from -4 to 15, y from -1 to 1.6
\put {$\bullet$} at -3 0
\put {$\bullet$} at -2 0
\put {$\bullet$} at -1 0
\put {$\bullet$} at 0 0
\put {$\bullet$} at 1 0
\put {$\bullet$} at 2 0

\put {$\bullet$} at -3 1
\put {$\bullet$} at -2 1
\put {$\bullet$} at -1 1
\put {$\bullet$} at 0 1
\put {$\bullet$} at 1 1
\put {$\bullet$} at 2 1

\put {$\bullet$} at -3 2
\put {$\bullet$} at -2 2
\put {$\bullet$} at -1 2
\put {$\bullet$} at 0 2
\put {$\bullet$} at 1 2
\put {$\bullet$} at 2 2

\setquadratic
\plot -3 2 -2.5 2.2 -2 2 /
\plot -1 2 -0.5 2.2 0 2 /
\plot 1 2 1.5 2.2 2 2 /
\plot -3 2 -2.5 1.8 -2 2 /
\plot -1 2 -0.5 1.8 0 2 /
\plot 1 2 1.5 1.8 2 2 /

\setlinear 
\plot -3 2 -3 -0.5 / 
\plot -2 2 -2 -0.5 / 
\plot -1 2 -1 -0.5 / 
\plot 0 2 0 -0.5  /
\plot -2 2 -1 2 /  
\plot 0 2 1 2 /
\plot 2 2 2.5 2 /
\plot -3 1 2.5 1 / 
\plot -3 0 2.5 0 / 
\plot 1 2 1 -0.5 /
\plot 2 2 2 -0.5 /

\put {$v_{1,1}$} [v] at -3.3 2 
\put {$v_{2,1}$} [v] at -3.3 1 
\put {$v_{3,1}$} [v] at -3.3 0 

\put {$v_{1,2}$} [v] at -1.9 2.1
\put {$v_{1,3}$} [v] at -1.2 2.1 
\put {$v_{1,4}$} [v] at 0.18 2.1
 
\put {$v_{2,2}$} [v] at -1.8 1.1 
\put {$v_{3,2}$} [v] at -1.8 0.1 
\put {$v_{2,3}$} [v] at -0.8 1.1

\put {$g_{1,1}$} [l] at -3.4 1.5
\put {$g_{2,1}$} [l] at -3.4 0.5
\put {$g_{2,2}$} [l] at -2.4 0.5
\put {$g_{2,3}$} [l] at -1.4 0.5

\put {$e_1$} [v] at -2.5 2.3
\put {$e_2$} [v] at -0.5 2.3
\put {$e_3$} [v] at 1.5 2.3

\put {$f_{1,2}$} [v] at -1.5 1.85
\put {$f_{1,4}$} [v] at 0.5 1.85
\put {$f_{1,1}$} [v] at -2.5 1.65
\put {$f_{1,3}$} [v] at -0.5 1.65
\put {$f_{1,5}$} [v] at 1.5 1.65

\put {$f_{2,1}$} [v] at -2.5 0.85
\put {$f_{3,1}$} [v] at -2.5 -0.2 
\put {$f_{3,2}$} [v] at -1.5 -0.2 

\put {$f_{2,2}$} [l] at -1.5 0.85
\put {$f_{2,3}$} [l] at -0.5 0.85

\arrow <0.235cm> [0.2,0.6] from -2.6 1 to -2.5 1 
\arrow <0.235cm> [0.2,0.6] from -1.6 1 to -1.5 1
\arrow <0.235cm> [0.2,0.6] from -0.6 1 to -0.5 1
\arrow <0.235cm> [0.2,0.6] from 0.4 1 to 0.5 1 
\arrow <0.235cm> [0.2,0.6] from 1.4 1 to 1.5 1

\arrow <0.235cm> [0.2,0.6] from -2.6 0 to -2.5 0 
\arrow <0.235cm> [0.2,0.6] from -1.6 0 to -1.5 0
\arrow <0.235cm> [0.2,0.6] from -0.6 0 to -0.5 0
\arrow <0.235cm> [0.2,0.6] from 0.4 0 to 0.5 0 
\arrow <0.235cm> [0.2,0.6] from 1.4 0 to 1.5 0

\arrow <0.235cm> [0.2,0.6] from -3 1.5 to -3 1.4 
\arrow <0.235cm> [0.2,0.6] from -3 0.5 to -3 0.4
\arrow <0.235cm> [0.2,0.6] from -3 -0.5 to -3 -0.6

\arrow <0.235cm> [0.2,0.6] from -2 1.5 to -2 1.4 
\arrow <0.235cm> [0.2,0.6] from -2 0.5 to -2 0.4
\arrow <0.235cm> [0.2,0.6] from -2 -0.5 to -2 -0.6

\arrow <0.235cm> [0.2,0.6] from -1 1.5 to -1 1.4 
\arrow <0.235cm> [0.2,0.6] from -1 0.5 to -1 0.4
\arrow <0.235cm> [0.2,0.6] from -1 -0.5 to -1 -0.6

\arrow <0.235cm> [0.2,0.6] from 0 1.5 to 0 1.4 
\arrow <0.235cm> [0.2,0.6] from 0 0.5 to 0 0.4
\arrow <0.235cm> [0.2,0.6] from 0 -0.5 to 0 -0.6

\arrow <0.235cm> [0.2,0.6] from 1 1.5 to 1 1.4 
\arrow <0.235cm> [0.2,0.6] from 1 0.5 to 1 0.4
\arrow <0.235cm> [0.2,0.6] from 1 -0.5 to 1 -0.6

\arrow <0.235cm> [0.2,0.6] from 2 1.5 to 2 1.4 
\arrow <0.235cm> [0.2,0.6] from 2 0.5 to 2 0.4
\arrow <0.235cm> [0.2,0.6] from 2 -0.5 to 2 -0.6

\arrow <0.235cm> [0.2,0.6] from 2.5 2 to 2.6 2 
\arrow <0.235cm> [0.2,0.6] from 2.5 1 to 2.6 1 
\arrow <0.235cm> [0.2,0.6] from 2.5 0 to 2.6 0

\arrow <0.235cm> [0.2,0.6] from -1.6 2 to -1.5 2 
\arrow <0.235cm> [0.2,0.6] from 0.4 2 to 0.5 2

\arrow <0.235cm> [0.2,0.6] from -2.7 2.175 to -2.8 2.12
\arrow <0.235cm> [0.2,0.6] from -2.4 1.8 to -2.3 1.84

\arrow <0.235cm> [0.2,0.6] from -0.7 2.175 to -0.8 2.12
\arrow <0.235cm> [0.2,0.6] from -0.4 1.8 to -0.3 1.84

\arrow <0.235cm> [0.2,0.6] from 1.3 2.175 to 1.2 2.12
\arrow <0.235cm> [0.2,0.6] from 1.6 1.8 to 1.7 1.84

\put {$\vdots$} at -3 -0.7
\put {$\vdots$} at -2 -0.7
\put {$\vdots$} at -1 -0.7
\put {$\vdots$} at 0 -0.7
\put {$\vdots$} at 1 -0.7
\put {$\vdots$} at 2 -0.7

\put {$\ldots$} at 2.9 2
\put {$\ldots$} at 2.9 1
\put {$\ldots$} at 2.9 0
\endpicture \] 
This is an infinite row-finite graph which does not satisfy 
Condition (K). There are four families of maximal tails
indexed by the integers $n\geq 1$: 
\begin{align*} 
M_n & =  \{v_{i,j}:1\leq i\leq n,\:1\leq j<\infty\}, \\ 
M^{2n-1} & =  \{v_{i,j}:1\leq i<\infty,\:
   1\leq j\leq 2n-1\}\cup\{v_{1,2n}\}, \\ 
M^{2n} & =  \{v_{i,j}:1\leq i<\infty,\:1\leq j\leq 2n\}, \\ 
R_n & =  \{v_{1,j}:1\leq j\leq 2n\}.  
\end{align*} 
In addition, $E^0$ is a maximal tail. Each maximal tail 
$R_n$ contains a loop without exits. 
On the other hand, all loops in 
$M_n$ and $M^n$ have exits. 
Since $E$ is row-finite there are no breaking vertices. 
Thus the gauge-invariant primitive ideals in $C^*(E)$ 
are $I_{\Omega(M_n)}$, $I_{\Omega(M^n)}$ and $\{0\}$. 
\end{example}

\begin{example}\label{example4} 
The following infinite graph $E$ is row-finite but does not 
satisfy Condition (K). 

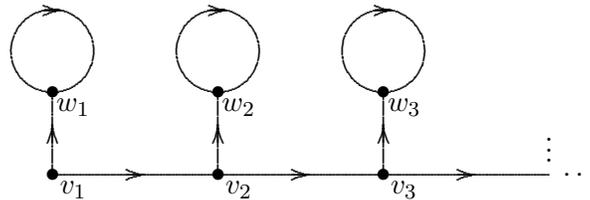
\[ \beginpicture
\setcoordinatesystem units <1.1cm,1.1cm>
\setplotarea x from -1 to 5.5, y from -0.25 to 1.9   

\put {$\bullet$} at 0 0
\put {$\bullet$} at 2 0
\put {$\bullet$} at 4 0
\put {$\bullet$} at 0 1
\put {$\bullet$} at 2 1
\put {$\bullet$} at 4 1

\setlinear 
\plot  0 1  0 0  4 0  4 1 / 
\plot  4 0  6 0 /
\plot  2 0  2 1 /  

\circulararc 360 degrees from 0 1 center at 0 1.5   
\circulararc 360 degrees from 2 1 center at 2 1.5   
\circulararc 360 degrees from 4 1 center at 4 1.5   

\arrow <0.235cm> [0.2,0.6] from 0.9 0 to 1.1 0 
\arrow <0.235cm> [0.2,0.6] from 2.9 0 to 3.1 0 
\arrow <0.235cm> [0.2,0.6] from 0 0.4 to 0 0.6  
\arrow <0.235cm> [0.2,0.6] from 2 0.4 to 2 0.6 
\arrow <0.235cm> [0.2,0.6] from 4.9 0 to 5.1 0
\arrow <0.235cm> [0.2,0.6] from 4 0.4 to 4 0.6  
\arrow <0.235cm> [0.2,0.6] from -0.1 1.99 to 0.1 1.99  
\arrow <0.235cm> [0.2,0.6] from 1.9 1.99 to 2.1 1.99  
\arrow <0.235cm> [0.2,0.6] from 3.9 1.99 to 4.1 1.99  

\put {$v_1$} at 0.25 -0.18 
\put {$v_2$} at 2.25 -0.18 
\put {$v_3$} at 4.25 -0.18 
\put {$w_1$} at 0.25 0.82 
\put {$w_2$} at 2.25 0.82 
\put {$w_3$} at 4.25 0.82 

\put {$\cdots$} at 6.4 0  
\put {$\vdots$} at 6 0.4  

\endpicture \] 
The maximal tails are $M=\{v_i:i\geq1\}$ and $M_n=\{w_n,v_1,\ldots,v_n\}$ 
for all $n\geq1$. Each $M_n$ contains a loop without exits, but $M$ 
does not contain any loops. Thus $I_{\Omega(M)}$ is the only
gauge-invariant primitive ideal  in $C^*(E)$. 
\end{example}

\begin{example}\label{example3} 
It was suggested in \cite[Remark 3.11]{fmr} that describing 
the ideals of the $C^*$-algebra of the following graph 
would be an interesting test question: 

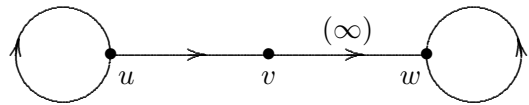
\[ \beginpicture
\setcoordinatesystem units <2.1cm,2.1cm>
\setplotarea x from -3 to 2, y from -0.4 to 0.1
\put {$\bullet$} at -2 0
\put {$\bullet$} at -1 0
\put {$\bullet$} at 0 0
\put {$(\infty)$} at -0.5 0.15

\setlinear 
\plot -2 0 0 0 / 
\arrow <0.235cm> [0.2,0.6] from -1.5 0 to -1.4 0 
\arrow <0.235cm> [0.2,0.6] from -0.5 0 to -0.4 0

\circulararc 360 degrees from -2 0 center at -2.3 0  
\circulararc 360 degrees from 0 0 center at 0.3 0 

\arrow <0.235cm> [0.2,0.6] from -2.6 0 to -2.58 0.1 
\arrow <0.235cm> [0.2,0.6] from 0.6 0 to 0.58 0.1

\put {$u$} [v] at -1.9  -0.15 
\put {$v$} [v] at -1 -0.15 
\put {$w$} [v] at -0.1 -0.15 
 
\endpicture \] 
This graph is not row-finite and does not satisfy Condition (K). There are no
breaking vertices, and there are three maximal  tails: $M_1=\{u\}$,
$M_2=\{u,v\}$ and $M_3=\{u,v,w\}$.  All loops in $M_2$ have exits, but 
$M_1$ and $M_3$ contain loops without exits. Thus there is exactly one 
gauge-invariant primitive ideal $I_w$, which corresponds to $M_2$, . 

We have $p_w C^*(E)p_w\cong C(\T)$ and there are infinitely many 
paths ending at $w$, so $I_w\cong C(\T)\otimes\K$. 
By Proposition~\ref{quotientalgebra}, the quotient $C^*(E)/I_w\cong 
C^*(E\setminus\{w\})$ is isomorphic to the Toeplitz algebra. 
\end{example}

\section{An application to $K$-theory}

One of our original motivations for analysing ideals in graph algebras 
was to extend \cite[Theorem~3.2]{rs} to arbitrary graphs, and we now do this;
the only difference between Theorem~\ref{Kgroups} below and
\cite[Theorem~3.2]{rs}  is the definition of $W$. Recall that
$\hat{\gamma}$ is the dual action of $\Z=\hat\T$ on the
crossed product $C^*(E)\rtimes_\gamma
\T$ by the gauge action $\gamma$, and that applying the integrated form of the 
canonical embedding $u:\T\rightarrow M(C^*(E)\rtimes_\gamma \T)$ to the
function $z\in L^1(\T)$  yields a projection $\chi_1=\int z u_z\,dz\in
C^*(E)\rtimes_\gamma
\T$. 

\begin{thm}\label{Kgroups}
Let $E$ be a graph, let $W=\{w\in E^0:
s^{-1}(w) \text{ is either empty or infinite}\}$, and let
$V=E^0\setminus W$. With respect to the decomposition $E^0=V\cup W$, the
$E^0\times E^0$ vertex matrix
\[
M(v,w):=\#\{e\in E^1: s(e)=v \;\text{ and }\; r(e)=w\},
\]
takes the block form
\[
M= \left( \begin{array}{ll} B & C \\ * & * \end{array} \right)
\]
where $B$ and $C$ have nonnegative integer entries. We define
$K:\Z^V\to\Z^V\oplus\Z^W$ by $K(x)=\big((1-B^t)x,-C^tx\big)$,
and $\phi:\Z^V\oplus\Z^W\to
K_0(C^*(E)\rtimes_\gamma\T)$ in terms of the usual basis by
$\phi(v)=[p_v\chi_1]$. Then $\phi$ restricts 
to an isomorphism $\phi|$ of
$\ker K$ onto $K_1(C^*(E))$, and induces an isomorphism
$\overline\phi$ of $\coker K$ onto $K_0(C^*(E))$
such that the following diagram commutes:
$$ \begin{diagram}\dgARROWLENGTH=0.8\dgARROWLENGTH
\node{\ker K} \arrow{e} \arrow{s,l}{\phi|}
\node{\Z^V}\arrow{e,t}{K} \arrow{s,l}{\phi}
\node{\Z^V\oplus\Z^{W}} \arrow{e} \arrow{s,l}{\phi}
\node{\coker K}\arrow{s,l}{\overline\phi}
\\
\node{K_1(C^*(E))} \arrow{e}
\node{K_0(C^*(E)\rtimes_\gamma \T)}\arrow{e,t}{1-\hat{\gamma}_*^{-1}} 
\node{K_0(C^*(E)\rtimes_\gamma \T)} \arrow{e}
\node{K_0(C^*(E)).}
\end{diagram} $$ 
\end{thm} 

Almost the entire proof of 
\cite[Theorem 3.2]{rs} applies in this
more general situation, and indeed only one point needs a different argument.
Recall that for integers
$m\leq n$ we denote by $E\times_1[m,n]$ the subgraph of $E\times_1\Z$
with vertices $\{(v,k):m\leq k\leq n,\;v\in E^0\}$
and edges $\{(e,k):m<k\leq n,\;e\in E^1\}$. It is 
essential for the proof of the theorem to know 
the $K$-theory of the corresponding algebra $C^*(E\times_1[m,n])$. 
We claim that $K_0(C^*(E\times_1[m,n]))$ is a free 
abelian group with free generators 
$$ \{[p_{(v,n)}]:v\in V\}\cup
   \{[p_{(v,k)}]:v\in W,\;m\leq k\leq n\}. $$ 
In the row-finite case this algebra is a direct sum of copies of the 
compacts (on Hilbert spaces of varying dimensions), and the claim 
is quite obvious. In general it is an $AF$-algebra with 
a more complicated structure. However, since any path 
in $E\times_1[m,n]$ has length at most $n-m$ the 
following lemma applies. 

\begin{lemma}\label{shortgraph} 
Let $E$ be a directed graph such that the length of any 
path $\alpha\in E^*$ does not exceed a fixed number $d$. Then 
$K_1(C^*(E))=0$ and $K_0(C^*(E))$ is the free abelian group generated by
$$ \{[p_v]:s^{-1}(v) \text{ is either empty or infinite}\}. $$ 
\end{lemma} 
\begin{proof} 
We proceed by induction on $d$. If $d=0$ then 
$C^*(E)$ is a direct sum of copies of $\C$, and the claim is clear. 
Suppose the lemma holds for all graphs with the maximum 
path length at most $d$, and let $E$ be a graph with 
the maximum path length at most $d+1$. We denote by $F$ 
the set of sinks in $E$. By Lemma~\ref{quotientalgebra} 
$C^*(E)/I_{\Sigma(F)}\cong C^*(E\setminus\Sigma(F))$, and thus 
there is an exact sequence 
\begin{equation} 
0\longrightarrow I_{\Sigma(F)}\longrightarrow C^*(E)
   \longrightarrow C^*(E\setminus\Sigma(F))\longrightarrow 0. 
\label{d-exact} 
\end{equation}  
In the graph $E\setminus\Sigma(F)$, the maximum path length is at most $d$.
Thus  the inductive hypothesis implies that $K_1(C^*(E\setminus
\Sigma(F)))=0$ and $K_0(C^*(E\setminus\Sigma(F)))$ 
is free abelian with free generators corresponding to 
vertices which emit infinitely many edges and sinks in 
$E\setminus\Sigma(F)$. Moreover, the sinks in this quotient graph 
are the vertices $v\in E$ which emit infinitely many edges 
but for which there is no path in $E$ from $v$ to another vertex which  
emits infinitely many edges. If $v\in F$ then the ideal of $C^*(E)$ 
generated by $p_v$ is isomorphic to $\K$, with 
$\{s_\mu s_\nu^*:\mu,\nu\in E^*,r(\mu)=r(\nu)=v$\} as a system 
of matrix units. Moreover, any two such ideals corresponding to 
different sinks have trivial intersection. Thus
$I_{\Sigma(F)}$ is the direct sum of these ideals,
$K_1(I_{\Sigma(F)})=0$, and $K_0(I_{\Sigma(F)})$ is the free abelian 
group generated by $\{[p_v]:v\in F\}$. Thus the 
six-term exact sequence of $K$-theory associated to (\ref{d-exact}) yields 
$K_1(C^*(E))=0$ and the short exact sequence 
\begin{equation}
0\longrightarrow K_0(I_{\Sigma(F)})\longrightarrow K_0(C^*(E))
\longrightarrow K_0(C^*(E\setminus\Sigma(F)))\longrightarrow 0. 
\label{splitKtheory}
\end{equation}  
Since $K_0(C^*(E\setminus\Sigma(F)))$ is free abelian, the sequence 
(\ref{splitKtheory}) splits. Furthermore, such a splitting 
$K_0(C^*(E\setminus\Sigma(F)))\rightarrow K_0(C^*(E))$ may be 
determined by lifting free generators; we choose to lift  $[p_v]$ to $[p_v]$.
This completes the proof of the  inductive step and the lemma. 
\end{proof} 

The rest of the proof of Theorem \ref{Kgroups} is exactly the same as 
in \cite[Theorem 3.2]{rs}. 

\begin{remark}
Of course this computation of $K$-theory is not entirely new, though we
believe the approach taken in \cite{rs} has much to commend it.  
Cuntz's original calculation of $K$-theory for Cuntz-Krieger algebras applies as
it stands to finite graphs without sinks or sources in which every loop  has an
exit
\cite[Proposition~3.1]{c}. This was extended to  locally finite graphs in
\cite{pr} and \cite{kprr}, and to arbitrary row-finite graphs in
\cite[Theorem~3.2]{rs};  for non-row-finite graphs without sinks or  sources,
we can apply the computations of
$K$-theory for the Cuntz-Krieger algebras of infinite matrices
(\cite[Theorem~4.5]{el2}, \cite[\S6]{sz0}, \cite[Theorem~4.1]{rs}).
Infinite graphs with finitely  many vertices are covered by \cite[Proposition
2]{sz2}, and arbitrary graphs  by \cite[Theorem~3.1]{dt2}, which was proved  by
reducing to the row-finite case and applying \cite[Theorem~3.2]{rs}.  

\end{remark}

\end{document}